\documentstyle[12pt]{article}

\newcommand{\const}{\mathop{\rm const}\limits}

\newcommand{\supp}{\mathop{\rm supp}\limits}

\begin{document}

\begin{center}

{\bf  AN IMBEDDING OF FRACTIONAL ORDER \\

\vspace{3mm}

SOBOLEV-GRAND LEBESGUE SPACES, } \\

\vspace{3mm}

{\bf with  constant evaluation.}\\

\vspace{3mm}

 $ {\bf E.Ostrovsky^a, \ \ L.Sirota^b } $ \\

\vspace{4mm}

$ ^a $ Corresponding Author. Department of Mathematics and computer science, Bar-Ilan University, 84105, Ramat Gan, Israel.\\

E - mail: \ eugostrovsky@list.ru\\

\vspace{3mm}

$ ^b $  Department of Mathematics and computer science. Bar-Ilan University,
84105, Ramat Gan, Israel.\\

\vspace{3mm}

E - mail:  sirota3@bezeqint.net \\

\end{center}

\vspace{3mm}

{\bf Abstract}. We extend in this article the classical imbedding theorems for fractional Lebesgue-Sobolev's
spaces into the so-called Grand Lebesgue spaces, with sharp constant evaluation. \par

\vspace{3mm}

\par

\vspace{3mm}

\par

\vspace{3mm}

\section{Introduction. Notations. Problem Statement.} \par

\vspace{3mm}

 For the (measurable) numerical function  $ f: R^n \to R $ the Fourier transform  $ F[f](\xi) = \hat{f}(\xi) $ is defined
as ordinary:

$$
 \hat{f}(\xi)  = F[f](\xi):=  \int_{R^n} e^{-i x\cdot \xi  } \ f(\xi) \ d x,
$$

 Hereafter $ x\cdot \xi $ denotes the inner (scalar) product of two vectors $ x, \ \xi \in R^d $ and $ |x| $ is ordinary Euclidean
norm of the vector $  x: \ |x| = \sqrt{x \cdot x}. $ \par

 The Fourier transform $ F[f](\xi) $  is correctly defined, e.g. if $ f \in \cup_{p \in [1,2]} L_p(R^n). $ The norm of the function
$ f $ in the Lebesgue, more exactly, Lebesgue-Riesz space $ L_p = L_p(R^n), \ p \ge 1 $  will be denoted for simplicity $ |f|_p: $

$$
|f|_p := \left[ \int_{R^n} |f(x)|^p \ dx  \right]^{1/p}.
$$
Let $ \Delta $ be the Laplacian.  The fractional, in general case, power $ \sqrt{- \Delta}^s $ may be defined as a pseudo-differential
operator through Fourier transform

$$
\sqrt{- \Delta}^s [f] := F^{-1}( |\cdot|^s \ F[f]).
$$

 The fractional Sobolev's space  $ W_p^s = W_p^s(R^n) $ consists by definition on all the functions $ f: R^n \to R   $
with finite norm (more precisely, semi-norm)

$$
||f||W_p^s \stackrel{def}{=} | (-\Delta)^{s/2} [f]|_p, \ p \ge 1,\eqno(1.1)
$$
the Aronszajn-Gagliardo norm; which is equivalent to the Slobodetskii $ ||\cdot||S_p^s  $  semi-norm:

$$
||f||S_p^s = ||f||S_p^s(R^n) \stackrel{def}{=} \left[ \int_{R^n} \int_{R^n} \frac{|f(x) - f(y)|^p}{|x-y|^{n + sp}} \ dx \ dy \right]^{1/p}. \eqno(1.2)
$$

 More information about these spaces, in particular on the imbedding theorem see in the works  \cite{Cotsiolis1}, \cite{Dyda1}, \cite{Frank1},
 \cite{Hurri1}, \cite{Ludwig1}, \cite{Maz'ya1}, \cite{Maz'ya2}, \cite{Nezza1}, \cite{Palatucci1}, \cite{Runst1}, \cite{Sloane1}, \cite{Slobodetskij1}, \cite{Uspenskii1}, \cite{Uspenskii2} etc. \par

\vspace{3mm}

{\bf Remark 1.1.} In the definition (1.2) instead the whole space $ R^n $ may be used arbitrary open set $ \Omega \subset R^n. $ In detail:

$$
 ||f||S_p^s(\Omega) \stackrel{def}{=} \left[ \int_{\Omega} \int_{\Omega} \frac{|f(x) - f(y)|^p}{|x-y|^{n + sp}} \ dx \ dy \right]^{1/p}. \eqno(1.2a)
$$

\vspace{3mm}

{\bf Remark 1.2.} In the case when $ s < 0 $ the fractional Laplace operator $ ( - \Delta)^{-|s|/2} $ coincides with Riesz potential

$$
( - \Delta)^{-|s|/2}[f](x) = I_{|s|}[f](x) = C_R(n,s) \int_{R^n} \frac{f(y)}{|x-y|^{n - |s|}} \ dy.
$$

 This case was detail investigated in many works, see e.g.\cite{Bogdan1}, \cite{Hardy1}, \cite{Lieb1}, \cite{Loss1},
\cite{Ostrovsky3}, \cite{Ostrovsky4}. Therefore, we can assume further $ s > 0. $\par

\vspace{3mm}

{\bf Remark 1.3.}  The "complete" norm in the fractional Sobolev's space may be introduced as follows.

$$
||f||V_p^s \stackrel{def}{=} \left[ |f|_p^p + | \ (-\Delta)^{s/2} [f] \ |_p^p \right]^{1/p} =
$$

$$
\left[ \int_{R^n} |f(x)|^p \ dx + (||f||W_p^s)^p \right]^{1/p}, \ p \ge 1. \eqno(1.3)
$$

\vspace{3mm}

{\bf Remark 1.4.}  The fractional Sobolev's  spaces are closely related with Besov spaces, see
\cite{Bennet1}, page 330-341. \par

\vspace{3mm}

 We will use further the following {\it sharp} imbedding theorem, see  \cite{Lieb1}, \cite{Lieb2}, section 4.3. Define the
following "constant"

$$
K(n,s) := \pi^{s/2} \frac{\Gamma((n-s)/2)}{\Gamma((n+s)/2)} \left\{ \frac{\Gamma(n)}{\Gamma(n/2)} \right\}^{s/n}, \eqno(1.4)
$$
 where $ \Gamma(\cdot) $ denotes the usually Gamma-function. If

 $$
   0 < s < n, \ 1 < p < n/s, \ u \in C_0^{\infty}(R^n), \   q = pn/(n - sp), \eqno(1.5)
 $$
 and $ q = pn/(n - sp), $  then

 $$
 |u|_q \le K(n,s) | \sqrt{ - \Delta }^s \ u |_p = K(n,s) ||u||W_p^s  \ - \eqno(1.6)
 $$
fractional Lebesgue-Sobolev's imbedding theorem.\par

 Note that the conditions (1.5) are also {\it necessary} for the inequality of a form (1.6) for some constant $ K(n,s). $
This assertion may be proved by means of the well-known scaling method, see e.g. \cite{Talenti1},   \cite{Ostrovsky28}.\par

  Evidently, the inequality (1.6) holds true for all functions $  u = u(x) $ belonging to the {\it completion} of the space
 $ C_0^{\infty} $ relative the fractional semi-norm $ ||\cdot|| W_p^s. $ \\

\vspace{3mm}

{\bf  Our aim in this  article is to extrapolate the fractional Lebesgue-Sobolev's imbedding theorem (1.6) into
the so-called fractional Grand Lebesgu-Sobolev's imbedding spaces,
and as a particular case the - into the so-called Exponential Orlicz  Spaces (EOS). }\\

\vspace{3mm}

A particular (but weight) case was considered in the previous  article \cite{Ostrovsky28}.\par

\vspace{4mm}

\section{Grand Lebesgue Spaces and Sobolev-Grand Lebesgue Spaces.}

\vspace{3mm}

 Now we will describe using Grand Lebesgue Spaces (GLS)  and Sobolev's  Grand Lebesgue Spaces (SGLS).\\

\vspace{3mm}

{\bf 1. Grand Lebesgue Spaces.}\par

\vspace{3mm}

 We recall in this section  for reader conventions some definitions and facts from the theory of GLS spaces.\par

\vspace{2mm}

 Recently, see \cite{Fiorenza1},  \cite{Fiorenza2}, \cite{Fiorenza3}, \cite{Iwaniec1}, \cite{Iwaniec2},
 \cite{Kozachenko1},\cite{Liflyand1}, \cite{Ostrovsky1}, \cite{Ostrovsky2},   etc.
 appears the so-called Grand Lebesgue Spaces $ GLS = G(\psi) =G\psi =
 G(\psi; A,B), \ A,B = \const, A \ge 1, A < B \le \infty, $ spaces consisting
 on all the measurable functions $ f: R^n \to R, $ (or more generally $ f: \Omega \to R) $ with finite norms

$$
     ||f||G(\psi) \stackrel{def}{=} \sup_{p \in (A,B)} \left[ |f|_p /\psi(p) \right].
     \eqno(2.1)
$$

Here $ \psi(\cdot) $ is some continuous positive on the {\it open} interval  $ (A,B) $ function such that

$$
     \inf_{p \in (A,B)} \psi(p) > 0, \ \psi(p) = \infty, \ p \notin (A,B).
$$
We will denote
$$
 \supp (\psi) \stackrel{def}{=} (A,B) = \{p: \psi(p) < \infty, \} \eqno(2.2)
$$

The set of all $ \psi $  functions with support $ \supp (\psi)= (A,B) $ will be denoted by $ \Psi(A,B). $ \par
  This spaces are rearrangement invariant, see \cite{Bennet1}, and
  are used, for example, in the theory of probability  \cite{Kozachenko1},
  \cite{Ostrovsky1}, \cite{Ostrovsky2}; theory of Partial Differential Equations \cite{Fiorenza2}, \cite{Iwaniec2};  functional analysis \cite{Fiorenza3},
  \cite{Iwaniec1},  \cite{Liflyand1}, \cite{Ostrovsky2}, \cite{Ostrovsky16}; theory of Fourier series \cite{Ostrovsky1}, theory of martingales \cite{Ostrovsky2},  mathematical statistics \cite{Sirota1}, \cite{Sirota2};  theory of approximation \cite{Ostrovsky7}   etc.\par
 Notice that in the case when $ \psi(\cdot) \in \Psi(A,B),  $ a function
 $ p \to p \cdot \log \psi(p) $ is convex, and  $ B = \infty, $ then the space
$ G\psi $ coincides with some {\it exponential} Orlicz space. \par
 Conversely, if $ B < \infty, $ then the space $ G\psi(A,B) $ does  not coincides with
 the classical rearrangement invariant spaces: Orlicz, Lorentz, Marcinkiewicz etc.\par

\vspace{3mm}

{\bf Remark 2.1.} If we define the {\it degenerate } $ \psi_r(p), r = \const \ge 1 $ function as follows:
$$
\psi_r(p) = \infty, \ p \ne r; \psi_r(r) = 1
$$
and agree $ C/\infty = 0, C = \const > 0, $ then the $ G\psi_r(\cdot) $ space coincides
with the classical Lebesgue space $ L_r. $ \par

\vspace{3mm}

{\bf Remark 2.2.} Let $ \xi: \Omega \to R $ be some (measurable) function from the set
$ L(p_1, p_2), \ 1 \le p_1 < p_2 \le \infty. $ We can introduce the so-called
{\it natural} choice $ \psi_{\xi}(p)$  as as follows:

$$
\psi_{\xi}(p) \stackrel{def}{=} |\xi|_p; \ p \in (p_1,p_2).
$$

\vspace{3mm}

{\bf 2. Sobolev-Grand Lebesgue  (SGL) spaces.}\par

\vspace{3mm}

 Let $ \psi = \psi(p) $ be the function described above.  We will say that the function $ f: R^d \to R $ belongs to the
 Sobolev-Grand Lebesgue space $ SGL\psi, $  iff the following it semi-norm is finite:

$$
||u||SGL\psi_s \stackrel{def}{=} \sup_{p \in \supp(\psi)} \left[ \frac{||u||W_p^s}{\psi(p)} \right]. \eqno(2.3)
$$

 This notion (up to equivalence)  for the integer values $ s $
appeared at first (presumably) in an article  \cite{Ostrovsky29} (2010); the near definition see in \cite{Onofrio1} (2013).\par

\vspace{3mm}

{\bf 3. Main result.} \\

\vspace{3mm}

{\bf Theorem 2.1.}  {\it Let } $ \psi(\cdot) \in \Psi(1,n/s) $ {\it where } $   0 < s < n, \ 1 < p < n/s. $
{\it  Define the function  }

$$
\nu(q) :=\psi\left( \frac{qn}{n+qs} \right), \eqno(2.4)
$$
{\it so that}

$$
\supp \nu(\cdot) = ( (n/(n-s), \ \infty  ).
$$

{\it Let also}

$$
 u \in C_0^{\infty}(R^n)\cap SGL\psi_s.  \eqno(2.4)
$$

{\it Proposition:}

$$
||u||G\nu \le K(n,s) \cdot ||u||SGL\psi_s, \eqno(2.5)
$$
{\it  where the constant $  K(n,s) $ is the best possible. }\par

\vspace{3mm}

{\bf Proof.} Let $ u \in C_0^{\infty}(R^n)\cap SGL\psi_s; $ we can and will suppose without loss of  generality
$ ||u||SGL\psi_s = 1. $ It follows from direct definition of Sobolev-Grand Lebesgue spaces

$$
||u||W_p^s \le \psi(p).
$$
 It follows from inequality (1.6)

$$
 |u|_q \le  K(n,s) ||u||W_p^s \le K(n,s) \psi(p).  \eqno(2.6)
$$

 Since  $  p = qn/(n+qs),$ we deduce from (2.6) for the values $ q > n/(n-s) $

$$
 |u|_q \le  K(n,s) \psi( qn/(n+qs))= K(n,s) \nu(q) = K(n,s) \ \nu(q) \ ||u||SGL\psi_s,\eqno(2.7)
$$
or equally

$$
||u||G\nu \le K(n,s) \ ||u||SGL\psi_s,
$$

 The {\it exactness } of the constant $ K(n,s) $ in (2.5) follows immediately from the one of main results, namely, theorem 2.1,
 of an article  \cite{Ostrovsky30}.\\

 Note that in the case when $ \psi(p) = \psi_r(p), \ r = \const \ge 1 $ we obtain as a particular case the ordinary fractional
 Sobolev's imbedding theorem.\par

\vspace{3mm}

\section{ Boundedness of fractional Laplacian in DGLS}

\vspace{3mm}

 Let us return to the fractional Sobolev's inequality (1.6):

$$
 |u|_q \le K(n,s) | \sqrt{ - \Delta }^s \ u |_p = K(n,s) ||u||W_p^s.  \eqno(1.6)
$$

 Assume now that the inequality (1.6) is true for some {\it interval  } of values $ s:  \ s \in (s_-, s_+), \ 0 < s_- < s_+ < 1.  $\par
More detail, let $  Q = Q(s_-, s_+)  $ be some (measurable) set  in the plane $ (p,s), \ p \ge 1, \ s \in (0,1), \ Q = \{ (p,s)  \} $ such that

$$
\forall s \in (s_-,s_+) \ \Rightarrow \exists p \ge 1, \ u \in W_p^s.  \eqno(3.1)
$$
 If for some values $ (p,s) \ u \notin W_p^s,  $ we denote formally $ ||u||W_p^s = \infty. $\par

 Denote

$$
Q_s = \{ p, \ p \ge 1, \ (p,s) \in Q  \};
$$
the "section" of the set $  Q  $  on the $ s \ - $ level.  Then $ \forall s \in (s_-, s_+) \ \Rightarrow Q_s \ne \emptyset. $\par
 The inequality (1.6) may be rewritten as follows:

$$
|u|_q \le K(n,s) \ ||u||W^s_{qn/(n+qs) }, \ s \in (s_-,s_+),
$$
therefore

$$
|u|_q \le  \inf_{s \in (s_-, s_+)} \left[ K(n,s) \ ||u||W^s_{qn/(n+qs) } \right]. \eqno(3.2)
$$

 The last inequality be reformulated on the language of Grand Lebesgue spaces as follows. Denote

$$
\zeta(q) = \inf_{s \in (s_-, s_+)} \left[ K(n,s) \ ||u||W^s_{qn/(n+qs) } \right], \ q \in (n/(n-s_+), \infty),
$$
then  $  ||u|| G\zeta \le 1. $ \par

\vspace{3mm}

{\bf  Definition of Derivative Grand Lebesgue spaces. } \\

\vspace{3mm}

 Let $ \tau = \tau(p,s), \ p > 1, \ s \in (s_-, s_+) $   be continuous function such that $  \inf_{p,s} \tau(p,s) = 1. $ By definition,
the function $ u = u(x), \ x \in R^n $ (or  $ x \in \Omega) $  belongs to the Derivative Grand Lebesgue  space $ DGL(\tau),$ if it has a
finite semi - norm

$$
|||u|||DGL\tau := \sup_{p > 1} \sup_{ s \in (s_-, s_+) } \left[ \frac{||u||W_p^s}{\tau(p,s)} \right]. \eqno(3.3)
$$

 We can now formulate the imbedding theorems in Derivative Grand Lebesgue spaces.\par

 \vspace{3mm}

{\bf Theorem 3.1.} {\it Let $  u(\cdot) \in DGL(\tau); $  then}

$$
|u|_q \le \inf_{s \in (s_-, s_+)} \left[ \ K(n,s) \ \tau(qn/( n+qs ), s) \ \right] \ \cdot |||u|||DGL\tau. \eqno(3.4)
$$

\vspace{3mm}

{\bf Proof} is alike one in the theorem 2.1. Indeed, let $ |||u|||DGL\tau = 1; $ then

$$
||u||W_p^s \le \tau(p,s), \ p = qn/(n + qs).
$$

 It remains to use the inequality (3.2) and take the minimum over $  s. $ \par

 \vspace{3mm}

 {\bf Corollary 3.1.} The inequality (3.4) may be reformulated on the language of GL spaces as follows. Denote

$$
 \lambda(q) = \inf_{s \in (s_-, s_+)} \left[ \ K(n,s) \ \tau(qn/( n+qs ), s) \ \right];
$$

then

$$
|u|_q \le \lambda(q) \cdot |||u|||DGL\tau,
$$
or equally

$$
 ||u||G\lambda \le |||u|||DGL\tau. \eqno(3.5)
$$

\vspace{3mm}

\section{ Weight generalization}

\vspace{3mm}

Let $ \Omega $ be an open {\it convex}  subset of a whole space $  R^n. $ Introduce  after
 R. L. Frank and R. Seiringer \cite{Frank1} (case   $ \Omega = R^n_+)  $ and M. Loss and C. Sloane
\cite{Loss1} (general case) the  following functions, measures  and operators:

$$
d_{\alpha}(x) :=\inf_{y \notin \Omega}  |x-y|^{\alpha},
$$

$$
D_{\alpha,n}(p) := 2 \pi^{(n-1)/2} \ \frac{\Gamma((1+\alpha)/2)}{\Gamma((n+\alpha)/2)} \
\int_0^1 \frac{\left|1-r^{(\alpha - 1)/p } \right|^p}{(1-r)^{1+\alpha} } \ dr, \ \alpha = \const \in (1,p);
$$

$$
g_{\alpha,n}(p) = \left[D_{\alpha,n}(p) \right]^{-1/p},
$$

$$
\mu_{\alpha}(A) := \int_A \frac{dx}{d_{\alpha}(x)}, \  A \subset \Omega;
$$

$$
\nu_{\alpha}(B) := \int \int_B \frac{dx \ dy}{|x-y|^{n + \alpha}}, \ B \subset \Omega \times \Omega.
$$

$$
\delta[f](x,y) := f(x) - f(y), \ f: R^n  \to R.\eqno(4.1)
$$

 The fractional weight Sobolev's type inequality

 $$
 |f|L_p(R^n,\mu_{\alpha}) \le g_{\alpha,n}(p) \ | \ \delta[f] \  |L_p(R^n \times R^n, \nu_{\alpha}) \eqno(4.2)
  $$
was proved by  R. L. Frank and R. Seiringer \cite{Frank1} (case   $ \Omega = R^n_+)  $ and M.Loss and C.Sloane
\cite{Loss1} (general case). See also \cite{Dyda1}. \par

\vspace{3mm}

{\bf Theorem 4.1.} {\it Let the function $ \delta[f](\cdot,\cdot)  $ belongs to some space $ G\psi $  on the set $ \Omega \times \Omega $
relative the measure $ \nu_{\alpha}. $ Put }

$$
\theta(p) = g_{\alpha,n}(p) \cdot \psi(p).
$$

{\it  Proposition: }

$$
||f||G\theta(\Omega, \mu_{\alpha}) \le 1 \cdot ||\delta[f]||G\psi(\Omega \times \Omega, \nu_{\alpha}), \eqno(4.3)
$$
{\it  where the constant "1" in (4.3) is the best possible. }\par

\vspace{3mm}

{\bf  Proof } is at the same as in theorem 2.1 and may be omitted. \par

\vspace{3mm}

\section{ Auxiliary results}

\vspace{3mm}

{\bf Constants $ K(n,s). $ }\\

 As long as  $ \Gamma(\epsilon) \sim 1/\epsilon, \ \epsilon \to 0+, $ we deduce at $ s \to n-0 $

$$
K(n,s) \sim \frac{2 \pi^{n/2}}{\Gamma(n/2)} \ \frac{1}{n-s}.
$$

 The point $ s = n- 0 $ is unique point of singularity for the function $ s \to K(n,s); $ for instance,
$ K(n, 0) = K(n,0+) = 1.$\\

\vspace{3mm}

{\bf Constants $ D_{\alpha,n}(p). $  }\\

 Denote for the values $ \alpha = \const > 1 $

$$
L_{\alpha}(p) = \int_0^1 \frac{\left|1-r^{(\alpha - 1)/p } \right|^p}{(1-r)^{1+\alpha} } \ dr, \ p \in (\alpha, \infty).
$$

Recall that

$$
D_{\alpha,n}(p) := 2 \pi^{(n-1)/2} \ \frac{\Gamma((1+\alpha)/2)}{\Gamma((n+\alpha)/2)} \ L_{\alpha}(p).
$$

The extreme points $ p = \alpha +0 $ and $  p\to \infty $ are points of singularity for the function $ p \to L_{\alpha}(p). $

\vspace{3mm}

{\bf A. Case } $ p \to \alpha + 0.  $ \\

 We have:

 $$
 r^{ ( \alpha - 1)/p} = e^{\ln r \cdot (\alpha - 1)/p }  \sim 1 -  \frac{|\ln r| (\alpha - 1)}{p},
 $$
therefore

$$
L_{\alpha}(p) \sim \frac{(\alpha - 1)^p}{p^p} \ \int_0^1 \frac{|\ln r|^p \ dr}{(1-r)^{1+\alpha}} \sim
$$

$$
\frac{(\alpha - 1)^p}{p^p} \ \int_0^1 (1-r)^{p-1-\alpha} \ dr =
\frac{(\alpha - 1)^p}{p^p} \cdot (p-\alpha)^{-1}.
$$

\vspace{3mm}

{\bf B. Case $ p \to \infty.  $ } \\

\vspace{3mm}

 We find:

$$
L_{\alpha}(p) \sim \frac{(\alpha - 1)^p}{p^p} \ \int_0^1 \frac{|\ln r|^p \ dr}{(1-r)^{1+\alpha}} \sim
$$

$$
\frac{(\alpha - 1)^p}{p^p} \int_0^1 |\ln r|^p \ dr =  \frac{(\alpha - 1)^p}{p^p} \cdot \Gamma(p+1).
$$

\vspace{3mm}

{\bf C. Non-asymptotical approach. } \\

\vspace{3mm}

 We will use the following elementary estimate:

$$
1 - \sinh(1) \ \epsilon \le e^{-\epsilon} \le  1 - \epsilon, \ \epsilon \in [0,1].
$$

 Let $ \Delta = \const \in (0,1);  $ for example, $ \Delta = \Delta_0 = 1/2. $ We calculate:

 $$
 J := \int_0^1 \frac{|\ln r|^p }{(1-r)^{1+\alpha}} dr  = \int_1^{\Delta} dr + \int_{\Delta}^1 dr = J_1 + J_2;
 $$

$$
J_1 \le (1-\Delta)^{-1 - \alpha} \int_0^{\Delta} |\ln r|^p \ dr \le
$$

$$
 (1-\Delta)^{-1 - \alpha} \int_0^1 |\ln r|^r \ dr  = \frac{\Gamma(p+1)}{( 1 - \Delta)^{1 + \alpha} )};
$$

$$
J_2 = \int_{\Delta}^1 \frac{|\ln r|^p }{(1-r)^{1+\alpha}} dr \le \frac{|\ln \Delta|^p}{( (1 - \Delta)^p  )} \ \frac{1}{p-\alpha},
$$
so

$$
J \le  \frac{\Gamma(p+1)}{( 1 - \Delta)^{1 + \alpha} } + \frac{|\ln \Delta|^p}{( (1 - \Delta)^p  )} \ \frac{1}{p-\alpha},
$$
following

$$
L_{\alpha}(p) \le \frac{(\alpha - 1)^p}{p^p} \cdot
\left[\frac{\Gamma(p+1)}{( 1 - \Delta)^{1 + \alpha} } + \frac{|\ln \Delta|^p}{( (1 - \Delta)^p  )} \ \frac{1}{p-\alpha} \right].
$$
 If we choose $ \Delta = 1/2, $ then

$$
L_{\alpha}(p) \le \frac{(\alpha - 1)^p}{p^p} \cdot
\left[ 2^{1 + \alpha} \Gamma(p+1) + (2 \ln 2)^p \ \frac{1}{p-\alpha} \right],
$$
and analogously

$$
L_{\alpha}(p) \ge C^p(\alpha) \cdot \frac{(\alpha - 1)^p}{p^p} \cdot
\left[ \Gamma(p+1) + \frac{1}{p-\alpha} \right], \ C(\alpha) \in (0,1).
$$

\vspace{3mm}

{\bf D. Constant of L.Cafarelli,  E. Valdinoci, O. Savin.  }\\

\vspace{3mm}

 There exists a "constant" $ Z = Z(n,s,p), \ s \in (0,1), \ p \ge 1 $ such that for all measurable set $  E \subset R^n $
with positive finite measure $ |E| $

$$
\int_{R^n \setminus E} \frac{dy}{|x-y|^{n + sp} } \ge Z (n,s,p) \ |E|^{-sp/n },
$$
see, e.g. \cite{Cafarelli1}, \cite{Savin1}. This constant play very important role in the theory of imbedding of
fractional Sobolev's spaces \cite{Nezza1}. \par

 We will understand as a capacity of the value $ Z (n,s,p) $ its maximal value, i.e.

 $$
 Z (n,s,p) \stackrel{def}{=} \inf_{x \in E} \
 \inf_{E: |E| \in (0,\infty)} \left\{ \int_{R^n \setminus E} \frac{dy}{|x-y|^{n + sp} } : |E|^{-sp/n}\right\}.
 $$

 Denote also as usually

$$
\omega_n = \frac{2 \pi^{n/2}}{\Gamma(n/2)} -
$$
the area of surface of unit sphere $  R^n;  $ recall that the volume of unit ball in this space is equal to $ \omega_n/n. $ \par

\vspace{3mm}

{\bf Proposition 5.1.}

$$
(sp)^{-1} \ \omega_n^{1 + sp/n } \ n^{-1 - sp/n  } \le  Z (n,s,p) \le  (sp)^{-1} \ \omega_n^{1 + sp/n } \ n^{- sp/n}.
$$

 The left - hand side follows immediately from lemma 6.1 in the article \cite{Nezza1} after simple calculations;  the right - hand side
may be obtained by choosing $  x = 0 $  and $ E = \{y: |y| \le 1. \} $ \par
 Obviously, the upper bound in the last inequality is attainable. \\

\vspace{3mm}

{\bf Acknowledgements.}

\vspace{3mm}

The authors are  grateful to  prof. S.V.Astashkin for useful conversations and to
Eleonora Di Nezza, Giampiero Palatucci, Enrico Valdinoci  for sending their remarkable articles.\\

\end{document}